\documentclass[conference]{IEEEtran}
\usepackage{stmaryrd}
\usepackage{amsfonts}
\usepackage{graphicx,times,psfig,amsmath,algorithmicx,algpseudocode,booktabs} % Add all your packages here

\IEEEoverridecommandlockouts    % to create the author's affliation portion
\begin{document}

% paper title: Must keep \ \\ \LARGE\bf in it to leave enough margin.
\title{\ \\ \LARGE\bf A new transformation into State Transition Algorithm for finding the global minimum\thanks{Xiaojun Zhou, Chunhua Yang and Weihua Gui are with the School of Information Science and Engineering, The Central South University, Changsha, Hunan, China, (email: tiezhongyu2010@gmail.com).} \thanks{The work was supported by the National Science Found for Distinguished Young Scholars of China (Grant No. 61025015).}}

\author{Xiaojun Zhou, Chunhua Yang and Weihua Gui}

% avoiding spaces at the end of the author lines is not a problem with
% conference papers because we don't use \thanks or \IEEEmembership
%\specialpapernotice{(Invited Paper)}

\maketitle

\begin{abstract}
To promote the global search ability of the original state transition algorithm, a new operator called axesion is suggested, which aims to search along the axes and strengthen single dimensional search. Several benchmark minimization problems are used to illustrate the advantages of the improved algorithm over other random search methods. The results of numerical experiments show that the new transformation can enhance the performance of the state transition algorithm and the new strategy is effective and reliable.
\end{abstract}

% no key words

\section{Introduction}
\PARstart{T}{he} basic random optimization(BRO) was proposed by Matyas\cite{cit:1} in 1965, and it is proved that the BRO can ensures convergence to a global minimum with probability one(\cite{cit:2},\cite{cit:3}). To  enhance the performance of the random optimization, various strategies have been introduced. In \cite{cit:4}, to adjust the parameters of mean and standard deviation in a Gaussian vector, a heuristic random optimization (HRO) war presented, utilizing two different mechanisms based on gradient information and reinforcement, respectively. In (\cite{cit:5},\cite{cit:6}), two approaches named adaptive random search technique(ARSET) and dynamic random search technique(DARSET) were put forward by Coskun Hamzacebi and Fevzi Kutay, to facilitate the determination of the global minimum. Then, another new random search algorithm named random selection walk(RSW) was raised through integrating the random selection and walk algorithms\cite{cit:7}.\\
\indent Based on the concepts of state and state transition, a metaheuristic random search method called state transition algorithm(STA) was emerged, and experimental results have shown that the proposed algorithm is advantageous for most of the test problems\cite{cit:8}. On the other hand, especially for optimization functions with independent variables, it is easy to get trapped into local optimum. In this study, a new operation call axesion transformation is raised to promote its search ability. Comparisons with other methods on several benchmark problems are presented in the paper, and the outcome is satisfactory, which shows that the proposed strategy is effective.
\section{the basic framework of STA}
In general, the form of state transition algorithm can be described as the following
\begin{equation}
\left \{ \begin{array}{ll}
x_{k+1}= A_{k}x_{k} + B_{k}u_{k}\\
y_{k+1}= f(x_{k+1})
\end{array} \right.
\end{equation}
where $x_{k}$ stands for a state, corresponding to a solution of the optimization problem; $A_{k}$ and
$B_{k}$ are state transition matrixes, which are usually transformation operators;
$u_{k}$ is the function with variables $x_{k}$ and history states; $f$ is the objective function or evaluation function.\\
\indent Using various types of space transformation for reference, three special
state transformation operators are defined to solve continuous function optimization problems.\\
(1) Rotation transformation
\begin{equation}
x_{k+1}=x_{k}+\alpha \frac{1}{n \|x_{k}\|_{2}} R_{r} x_{k}
\end{equation}
where $x_{k}$ $\in$ $\Re^{n\times1}$, $\alpha$ is a positive constant, called rotation factor;
$R_{r}$ $\in$ $\Re^{n\times n}$,is random matrix with its elements belonging to the range of [-1, 1]
and $\|\cdot\|_{2}$ is 2-norm of a vector. It has proved that the rotation transformation
has the function of searching in a hypersphere\cite{cit:8}.\\
\noindent (2) Translation transformation\\
\begin{equation}
x_{k+1} = x_{k}+  \beta  R_{t}  \frac{x_{k}-x_{k-1}}{\|x_{k}-x_{k-1}\|_{2}}
\end{equation}
where $\beta$ is a positive constant, called translation factor; $R_{t}$ $\in \Re^{1}$ is a random variable
with its elements belonging to the range of [0,1]. It has illustrated
the translation transformation has the function of searching along a line from $x_{k-1}$ to $x_{k}$
 at the starting point $x_{k}$, with the maximum length of $\beta$\cite{cit:8}.\\
(3) Expansion transformation
\begin{equation}
x_{k+1} = x_{k}+  \gamma  R_{e}x_{k}
\end{equation}
where $\gamma$ is a positive constant, called expansion factor; $R_{e} \in \Re^{n \times n}$is a random diagonal
matrix with its elements obeying the Gaussian distribution. It has also stated the expansion transformation
has the function of expanding the elements in $x_{k}$ to the range of [-$\infty$, +$\infty$],
searching in the whole space\cite{cit:8}.\\
\indent The procedures of the original state transition algorithm can be outlined in the following pseudocode.\\
\begin{algorithmic}[1]
\State Initialize feasible solution $x(0)$ randomly, set $\alpha, \beta, \gamma$, and $k \gets 0$
\Repeat
   \State $k \gets k +1 $
   \While{$\alpha \le tolerance$}
    \State $State \gets op\_rotate(x(k-1),SE,\alpha)$
    \If{min $f(State) < f(x(k-1))$}
        \State Updating $x(k-1)$
        \State $State \gets op\_translate(x(k-1),SE,\beta)$
        \If{min $f(State) < f(x(k-1))$}
            \State Updating $x(k-1)$
        \EndIf
    \EndIf
    \State $\alpha \gets \frac{\alpha}{\textit{fc}} $
   \EndWhile
   \State $State \gets op\_expand(x(k-1),SE,\gamma)$
   \If{min $f(State) < f(x(k-1))$}
    \State Updating $x(k-1)$
    \State $State \gets op\_translate(x(k-1),SE,\beta)$
    \If{min $f(State) < f(x(k-1))$}
        \State Updating $x(k-1)$
    \EndIf
   \EndIf
  \State $x(k) \gets x(k-1)$
\Until{the specified termination criterion is met}
\end{algorithmic}
\quad \\
where (\textit{SE}) is search enforcement, which means the times of the transformation. Operators such as $op\_rotate(\cdot), op\_translate(\cdot)$ and $op\_expand(\cdot)$ correspond
to the rotation, translation, and expansion, respectively. \textit{fc} is a constant coefficient used for lessening the $\alpha$.
By the way, the translation operator will only be performed when a better solution is obtained.

\section{A new transformation into the STA}
\begin{figure}[h!]
\centering
\includegraphics[width=6cm,height=6cm]{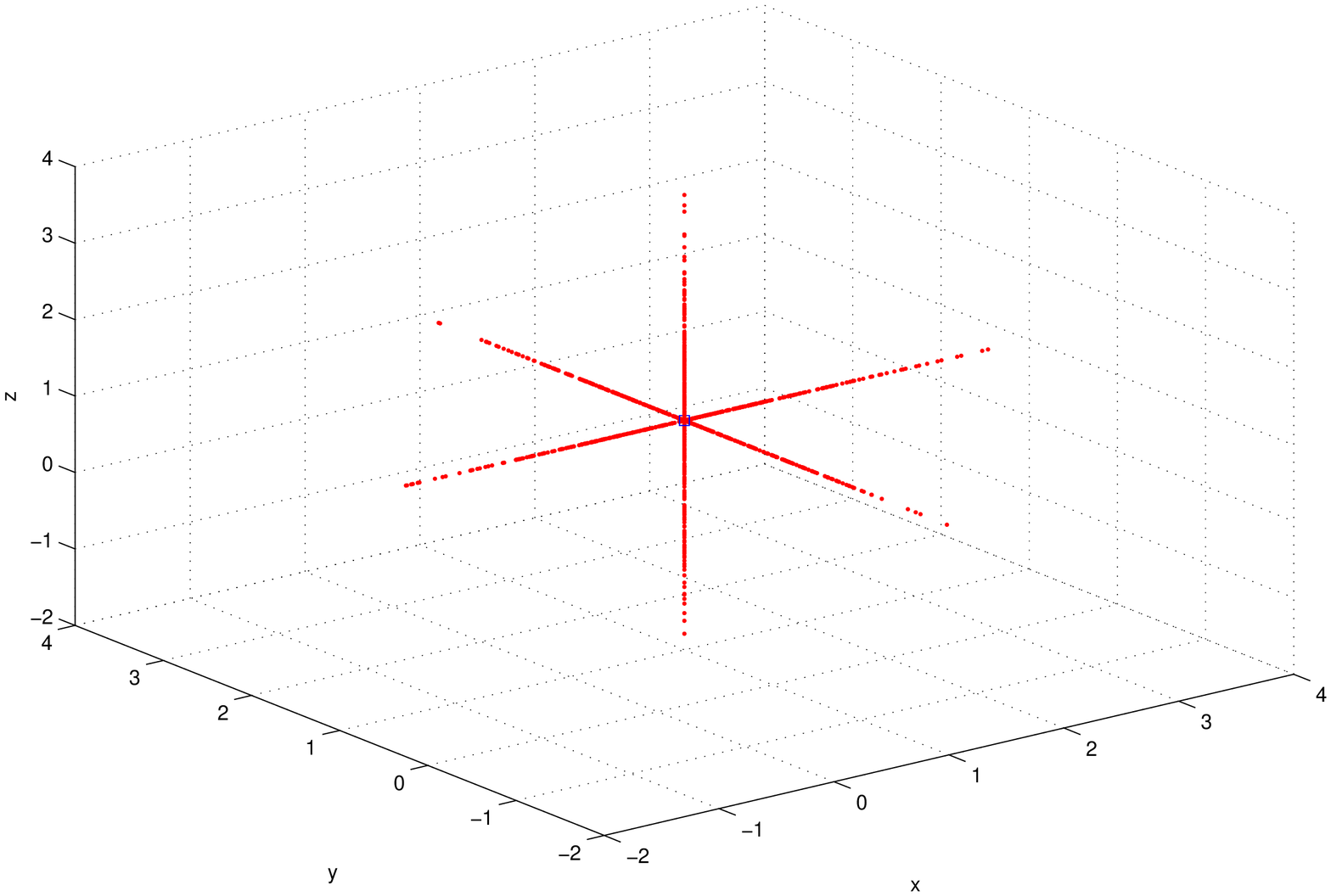}
\caption{the function of axesion transformation}
\label{tab:axesion}
\end{figure}

In the proposed STA, three different state transformation operators are designed, aiming for exploration(global search) and
exploitation(local search) as well as the equilibrium between them. It is beneficial for optimization functions with relevant variables;
however, for functions with independent variables, it is necessary to intensify the single dimension search.\\
\indent To simplify the one dimensional search, a new operator called axesion is added to the STA, which, in its meaning, aims to search along
each axes.\\
(4) Axesion transformation
\begin{equation}
x_{k+1} = x_{k}+  \delta  R_{a}  x_{k}\\
\end{equation}
where $\delta$ is a positive constant, called axesion factor; $R_{a}$ $\in \Re^{n \times n}$ is a random diagonal matrix with its
elements obeying the Gaussian distribution and only one random index has value. For example, $R_{a}$ $\in \Re^{3 \times 3}$ can be the following
styles:\\
\begin{center}
$\left(
  \begin{array}{ccc}
    0.1 & 0 & 0 \\
    0   & 0 & 0 \\
    0   & 0 & 0 \\
  \end{array}
\right)$
$\left(
  \begin{array}{ccc}
    0 & 0   & 0 \\
    0 & 0.2 & 0 \\
    0 & 0   & 0 \\
  \end{array}
\right)$
$\left(
  \begin{array}{ccc}
    0 & 0 & 0   \\
    0 & 0 & 0   \\
    0 & 0 & 0.3 \\
  \end{array}
\right)$
\end{center}
\quad \\
\indent To illustrate the functions of the axesion transformation, let suppose the $x_{k} = [1;1;1]$ and $\delta = 1$; then, after 1e3 times of independent axesion transformation, the distribution of $x_{k+1}$ can be illustrated on the space as described in Figure.{\ref{tab:axesion}}.\\
\indent When the new transformation is introduced into the original STA, it will follow the same procedures as rotation and expansion, while translation transformation will only be performed when a better solution is gained by axesion. In the meanwhile, to reduce the computational complexity, the rotation
transformation will be executed in an outer loop instead of an inner loop, that is to say, the $\alpha$ factor will vary in a periodic way.

\section{Numerical experiments}
Experiments are divided into two groups. The first group consist of five benchmark problems which have been tested by HRO, ARSET, and RSW, while the second group benchmark problems have been tested by DARSET and RSW. In the proposed original STA and new STA, all of the problems will be carried out. Parameters of the STAs are described in Table \ref{tab:STAs}, and the lessening coefficient \textit{fc} will be 4 in original STA and 2 in new STA.
\begin{table}[h!]
\begin{center}
\caption{Parameters of the STAs}
\label{tab:STAs}
\footnotesize
\begin{tabular}{{cc}}
\hline
\toprule[1pt]
\textit{Parameter} & \textit{Value}\\
\hline
epoch number                  & 1000 \\
$SE$(Search enforcement)      & 32   \\
$\alpha$                      & 1 $\rightarrow$ 1e-4  \\
$\beta$                       & 1    \\
$\gamma$                      & 1    \\
$\delta$                      & 1    \\
\textit{fc}                   & 4(2) \\
\bottomrule[1pt]
\hline
\end{tabular}
\end{center}
\end{table}
\subsection{The first Group}
\begin{figure}[h!]
\centering
\includegraphics[width=6cm,height=6cm]{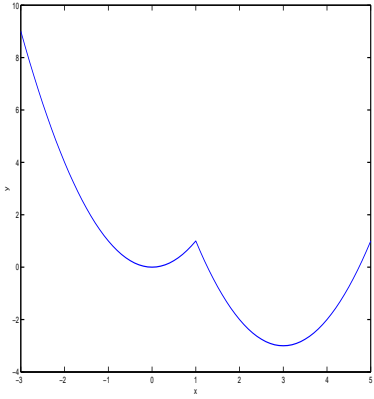}
\caption{the landscape of $f_1$}
\label{tab:f1}
\end{figure}

\begin{figure}[h!]
\centering
\includegraphics[width=6cm,height=6cm]{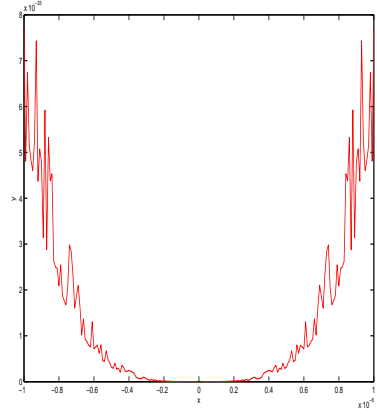}
\caption{the landscape of $f_2$}
\label{tab:f2}
\end{figure}

\begin{figure}[h!]
\centering
\includegraphics[width=6cm,height=6cm]{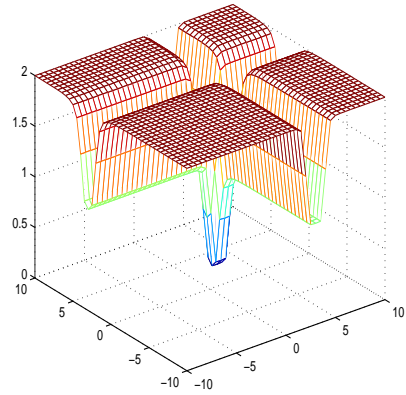}
\caption{the landscape of $f_3$}
\label{tab:f3}
\end{figure}

\begin{figure}[h!]
\centering
\includegraphics[width=6cm,height=6cm]{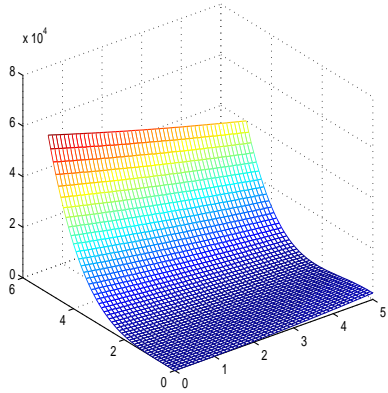}
\caption{the landscape of $f_4$}
\label{tab:f4}
\end{figure}

\begin{figure}[h!]
\centering
\includegraphics[width=6cm,height=6cm]{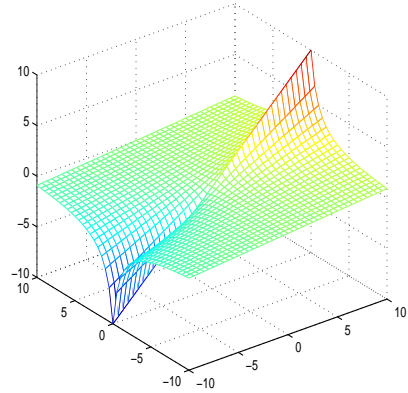}
\caption{the landscape of $f_5$}
\label{tab:f5}
\end{figure}

\noindent(1) Problem $f_1$\\
\indent The first problem is taken from [4,5,7], objective function of the problem is given as follows
\begin{equation}
f_1 = \left\{\begin{array}{ll}
x^2, & \textrm{if}\quad  x \leq 1, \\
(x-3)^2 -3, & \textrm{if}\quad x > 1.
\end{array}\right.
\label{tab:Eqf1}
\end{equation}
\indent As can be seen from Fig. \ref{tab:f1}, the function has two minimums, one lying on $x = 0$ and the other one on $x = 3$.
Their results of the HRO, ARSET, RSW, and STAs are given in Table \ref{tab:rf1}. Both STAs can achieve the global minimum in the end,
as a matter of fact, STAs can meet the optimum in no more than 10 epoches.
\quad
\begin{table}[h!]
\begin{center}
\caption{Results for the $f_1$ }
\label{tab:rf1}
\begin{tabular}{{cccc}}
\hline
\toprule[1pt]
\textit{Algorithms} & Best $x$ & Best $f(x)$\\
\hline
HRO           & 3.000324 &  -3 \\
ARSET         & 3        &  -3 \\
RSW           & 3        &  -3 \\
STA(original) & 3        &  -3 \\
STA(new)      & 3        &  -3 \\
\bottomrule[1pt]
\hline
\end{tabular}
\end{center}
\end{table}\\
\noindent(2) Problem $f_2$\\
\indent The second problem is taken from [4,5,7], objective function of the problem is given as follows
\begin{equation}
f_2 = [x sin(\frac{1}{x})]^4 + [x cos(\frac{1}{x})]^4, f(0)= \lim_{x\rightarrow 0 }f(x)=0
\label{tab:Eqf2}
\end{equation}
\indent As can be seen from Fig. \ref{tab:f2}, the function has numerous local minimums.
Their results of the HRO, ARSET, RSW, and STAs are given in Table \ref{tab:rf2}. Both STAs can achieve the global minimum in the end,
as a matter of fact, STAs can meet the optimum in no more than 50 epoches.
\begin{table}[h!]
\begin{center}
\caption{Results for the $f_2$ }
\label{tab:rf2}
\begin{tabular}{{cccc}}
\hline
\toprule[1pt]
\textit{Algorithms} & Best $x$ & Best $f(x)$\\
\hline
HRO           & 2.4000e-005 &  2.8595e-019 \\
ARSET         & -2.53e-011  &  2.21e-043   \\
RSW           & 8.17e-82    &  0           \\
STA(original) & 2.0447e-082 &  0           \\
STA(new)      & 3.5197e-084 &  0           \\
\bottomrule[1pt]
\hline
\end{tabular}
\end{center}
\end{table}\\
\noindent(3) Problem $f_3$\\
\indent The third problem is taken from [5,7], objective function of the problem is given as follows
\begin{equation}
f_3 = \frac{(x-3)^8}{1+(x-3)^8}+\frac{(y-3)^4}{1+(y-3)^4}
\label{tab:Eqf3}
\end{equation}
\indent As can be seen from Fig. \ref{tab:f3}, the function also has numerous local minimums.
Their results of the ARSET, RSW, and STAs are given in Table \ref{tab:rf3}. Compared with ARSET and RSW, STAs can get better
solution than them, in the same time, STAs can meet the specified precision in no more than 100 epoches.
\begin{table}[htbp]
\begin{center}
\caption{Results for the $f_3$ }
\label{tab:rf3}
\begin{tabular}{{cccc}}
\hline
\toprule[1pt]
\textit{Algorithms} & Best $x$ & Best $y$ & Best $f(x,y)$\\
\hline
ARSET         & 3.0015      &  3      &  5.04e-023  \\
RSW           & 2.9996      &  3      &  3.43e-28   \\
STA(original) & 3.0000      &  3.0000 &  5.8715e-033\\
STA(new)      & 3.0000      &  3.0000 &  1.0335e-035\\
\bottomrule[1pt]
\hline
\end{tabular}
\end{center}
\end{table}\\
\noindent(4) Problem $f_4$\\
\indent The fourth problem is taken from [5,7], objective function of the problem is given as follows
\begin{equation}
f_4 = 100(x-y^2)^2 + (1-x)^2
\label{tab:Eqf4}
\end{equation}
\indent Fig. \ref{tab:f4} shows the graph of the function, which is widely used for testing because there is a valley
in the landscape, making it hard to optimize. Their results of the ARSET, RSW, and STAs are given in Table \ref{tab:rf4}.
Compared with ARSET and RSW, STAs seem a little deficient for the function.
\begin{table}[h!]
\begin{center}
\caption{Results for the $f_4$ }
\label{tab:rf4}
\begin{tabular}{{cccc}}
\hline
\toprule[1pt]
\textit{Algorithms} & Best $x$ & Best $y$ & Best $f(x,y)$\\
\hline
ARSET         & 1      &  1      &  4.02e-016  \\
RSW           & 1      &  1      &  1.97e-31   \\
STA(original) & 1.0000 &  1.0000 &  8.2040e-012\\
STA(new)      & 1.0000 &  1.0000 &  3.7678e-012\\
\bottomrule[1pt]
\hline
\end{tabular}
\end{center}
\end{table}\\
\noindent(5) Problem $f_5$\\
\indent The fifth problem is taken from [5,7],objective function of the problem is given as follows
\begin{equation}
f_5 = \frac{x}{1+|y|}
\label{tab:Eqf5}
\end{equation}
\indent Fig. \ref{tab:f5} shows the graph of the function, which is indifferentiable at the minimum point. Their results of the ARSET, RSW, and STAs are given in Table \ref{tab:rf5}.
Compared with ARSET and RSW, STAs are much better than them, because only they can meet the global minimum.
\begin{table}[h!]
\begin{center}
\caption{Results for the $f_5$ }
\label{tab:rf5}
\begin{tabular}{{cccc}}
\hline
\toprule[1pt]
\textit{Algorithms} & Best $x$ & Best $y$ & Best $f(x,y)$\\
\hline
ARSET         & -10     &  6.67e-008 & -10    \\
RSW           & -9.9996 &  -6.57e-17 & -9.9996\\
STA(original) & -10.0000&  0.0000    & -10    \\
STA(new)      & -10.0000&  0.0000    & -10    \\
\bottomrule[1pt]
\hline
\end{tabular}
\end{center}
\end{table}
\subsection{The second Group}
Problems $g_1$ to $g_{15}$ are taken from [6] and [7], which have been tested by DARSET and RSW. The detailed information of
the test functions are listed in  Table \ref{tab:grp2}. In DARSET and RSW, epoch number
varies from 250200 to 2508000; however, in original STA and new STA, both the epoch number is fixed at 1000, and 10 independent times are run.\\
\indent As described in Table \ref{tab:rgrp}, the original STA can get the same results as DARSET and RSW in
$g_4, g_5,g_6,g_7,g_8,g_{14}$ and $g_{15}$. For $g_1$ and $g_{10}$, the results of original STA are better than DARSET but a little inferior to BSW, for $g_2$ and $g_3$, the original STA can get better results than BSW but a little inferior to DARSET, and for $g_{11}$, the original STA can gain better solution than both DARSET and BSW.\\
\indent On the other hand, the new STA can achieve the best results of all functions except for $g_9,g_{12}$ and $g_{13}$. Even so, it is obvious to find that the new STA actually can achieve the global minimums in a specified precision.\\
\indent Through the numerical experiments, it has proved that the new STA has better performance than the original STA in terms of the global search ability and the proposed axesion transformation is beneficial for STA.\\
\indent The reasons of the efficiency of the new strategy can be explained as: (1) axesion transformation enlarge the search space. Not only rotation but also expansion, both of them search in the neighborhood of the best solution, of which, one search in a unit hypersphere, the other search in a relative broader space but with low probability, while the new operator can also search in global space, that is to say, the probability of finding the global optimum is increased. (2) the new transformation has the function of single dimensional search, which is advantageous for functions with independent variables. (3) the new transformation can achieve good quality solution with high precision. When searching in a single direction, the axesion can enhance the depth of search.
\begin{table*}[htbp!]
\caption{The second group functions for test in this paper}
\label{tab:grp2}
\small %\tiny \footnotesize
\begin{tabular}{{lccc}}
\hline
\toprule[1pt]
Function & Dimension  & Variable range & Theoretical best \\
\hline
$g_1=x^2 + 2y^2 -0.3cos(3\pi x)-0.4cos(4\pi y)+0.7$                         & 2 & [-1.28,1.28]    & 0       \\
$g_2=[cos(2\pi x )+ cos(2.5\pi x )-2.1][2.1 - cos(3\pi y )+ cos(3.5\pi y)]$ & 2 & [-1,1]          & -16.0917\\
$g_3=[0.002+\sum_{j=1}^{25}(j+\sum_{i=1}^2(x_{i}-a_{ij})^6)^{-1}]^{-1}$     & 2 & [-65.536,65.536]& 0.9980  \\
$a = \left|\begin{array}{ccccccccccc}
 -32 &  -16 & 0 &  16 & 32 & ...&  -32 & -16 & 0 & 16 & 32\\
 -32 &  -32 &-32& -32 &-32 & ...& 32  & 32  & 3 & 32 & 32
\end{array}\right|$ & & & \\
$g_4=(y - \frac{5.1}{4 \pi^2}x^2 + \frac{5}{\pi}x -6)^2 + 10(1-\frac{1}{8\pi})cos(x)+10$ & 2
&
$
\begin{array}{c}
x \in [-5,10]\\
y \in [0,15]
\end{array}
$
& 0.3979\\
$g_5=(4-2.1x^2+\frac{x^4}{3})x^2 + xy + (4y^2-4)y^2$ & 2
&
$
\begin{array}{c}
x \in [-3,3]\\
y \in [-2,2]
\end{array}
$
& -1.0316\\
$
\begin{array}{l}
g_6=[1+(x+y+1)^2(19-14x+3x^2-14y+6xy+3y^2)] \\
\times [30+(2x-3y)^2(18-32x+12x^2+48y-36xy+27y^2)]
\end{array}
$
& 2 & [-5,5] & 3\\
$g_7=[\sum_{i=1}^5icos((i+1)x+i)]\times[\sum_{i=1}^5icos((i+1)y+i)]$                        & 2 & [-10,10]            & -186.7309 \\
$g_8=\sum_{i=1}^5(x(a_{i})^y(b_{i})^z-c_{i})^2$                                             & 3 & $[-\infty,\infty]$  & 8.0128    \\
$a=|5~ 3 ~0.6 ~0.1~ 3|,b=|10~1~0.6~2~1.8|,c=|2.122~9.429~23.57~74.25~6.286|$ & & & \\
$
\begin{array}{l}
g_9=100(x_2-x_1^2)^2+(1-x_1)^2+90(x_4-x_3^2)^2+(1-x_3)^2\\
+10.1[(x_2-1)^2+(x_4-1)^2]+19.8(x_2-1)(x_4-1)
\end{array}
$
& 4 & [-10,10]            & 0\\
$g_{10}=\sum_{i=1}^{19}[(x_i^2)^{(x_{i+1}^2+1)}+(x_{i+1}^2)^{(x_{i}^2+1)}]$                                               & 20 & [-1,4]   & 0   \\
$g_{11}=(\frac{\pi}{20})[10sin^2(\pi x_1)+\sum_{i=1}^{19}((x_i-1)^2(1+10sin^2(\pi x_{i+1}))+(x_{20}-1)^2]$                & 20 & [-10,10] & 0   \\
$g_{12}=100(y-x^2)^2+(1-x)^2$                                                                                             & 2  & [-10,10] & 0   \\
$g_{13}=exp(0.5(x^2+y^2-25)^2)+sin^4(4x-3y)+0.5(2x+y-10)^2$                                                               & 2  & [-5,5]   & 1   \\
$g_{14}=0.1[12+x^2+\frac{1+y^2}{x^2}+\frac{x^2y^2+100}{(xy)^4}]$                                                          & 2  & [0,10]   & 1.74\\
$g_{15}=(x_1+10x_2)^2+5(x_3-x_4)^2+(x_2-2x_3)^4+10(x_1-x_4)^4$                                                            & 4  & [-5,5]   & 0   \\
\bottomrule[1pt]
\hline
\end{tabular}
\end{table*}
\begin{table*}[htbp!]
\caption{Test results of the second group functions in this paper}
\label{tab:rgrp}
\small %\tiny \footnotesize
\begin{tabular}{{ccccccccc}}
\hline
\toprule[1pt]
Function & \multicolumn{2}{c}{DRASET} & \multicolumn{2}{c}{RSW} & \multicolumn{2}{c}{STA(original)} & \multicolumn{2}{c}{STA(new)} \\

              & Best & Average  & Best & Average & Best & Average  &Best & Average  \\
\hline
$g_1$  & 0        & 9.10e-016 & 0         & 0        & 0          & 5.3147e-012 & 0          & 0          \\
$g_2$  & -16.0917 & -16.0917  & -16.0917  & -15.7399 & -16.0917   & -16.0917    & -16.0917   & -16.0917   \\
$g_3$  & 0.998    & 1.5885    & 0.998     & 6.3728   & 0.9980     & 3.9354      & 0.9980     & 0.9980     \\
$g_4$  & 0.3979   & 0.3979    & 0.3979    & 0.3979   & 0.3979     & 0.3979      & 0.3979     & 0.3979     \\
$g_5$  & -1.0316  & -1.0316   & -1.0316   & -1.0316  & -1.0316    & -1.0316     & -1.0316    & -1.0316    \\
$g_6$  & 3        & 3         & 3         & 3        & 3.0000     & 3.0000      & 3.0000     & 3.0000     \\
$g_7$  & -186.7309& -186.7309 & -186.7309 & -186.7309& -186.7309  & -186.7309   & -186.7309  & -186.7309  \\
$g_8$  & 8.0128   & 8.0128    & 8.0128    & 8.0128   & 8.0128     & 8.0128      & 8.0128     & 8.0128     \\
$g_9$  & 3.72e-12 & 9.30e-06  & 1.28e-28  & 2.15e-28 & 2.8718e-010& 1.1802e-009 & 8.3086e-011& 1.1344e-009\\
$g_{10}$ & 2.45e-16 & 4.02e-15  & 0         & 0        & 0          & 0           & 4.9783e-094& 2.7247e-084\\
$g_{11}$ & 5.93e-12 & 26.227    & 2.36e-32  & 3.3927   & 7.2021e-011& 1.0417      & 2.6223e-011& 3.8022e-011\\
$g_{12}$ & 3.91e-15 & 4.28e-14  & 2.84e-29  & 6.07e-28 & 8.9683e-014& 3.8771e-012 & 9.5239e-014& 9.9002e-012\\
$g_{13}$ & 1        & 1.0077    & 1.0091    & 1.0091   & 1.0000     & 1.0375      & 1.0000     & 1.0225     \\
$g_{14}$ & 1.7442   & 1.7442    & 1.7442    & 1.7442   & 1.7442     & 1.7442      & 1.7442     & 1.7442     \\
$g_{15}$ & 8.17e-09 & 1.68e-07  & 1.02e-11  & 1.71e-11 & 2.1942e-014& 6.4995e-009 & 9.9870e-014& 1.0542e-007\\
\bottomrule[1pt]
\hline
\end{tabular}
\end{table*}
\section{Conclusions}
As a random search method, the original STA has shown the great ability in optimizing continuous functions. To enhance the global search capability of the STA, axesion transformation is introduced into the original STA, which aims to search along a single dimension in depth. The results of the numerical experiments have testified the efficiency and reliability of the new STA. Comparisons with other random optimization methods, the outcome of the experiments has also revealed the advantages of the STAs. By the way, other functions with independent variables have also been tested, and results are more satisfactory.



% http://www.ctan.org/tex-archive/macros/latex/contrib/supported/IEEEtran/testflow/bibtex
%\bibliographystyle{IEEEtran.bst}

\begin{thebibliography}{1}
\bibitem{cit:1} Matyas, J, ``Random optimization," {\it Automation and Remote Control}, vol. 26, pp. 246-253, 1965.
\bibitem{cit:2} Francisco J. Solis and Roger J-B. Wets, ``Minimization by random search techniques," {\it Mathematics of operations research}, vol. 6, no. 2, pp. 19-30, 1981.
\bibitem{cit:3} N.Baba, T.Shoman,Y.Sawaragi, ``A modified convergence therorem for a random optimization method," {\it Information Science}, vol. 13, pp. 159-166, 1977.
\bibitem{cit:4} Junyi Li and R. Russell Rhinehart, ``Heuristic random optimization," {\it Computers and Chemical Engineering}, vol. 22, no. 3, pp. 427-444, 1998.
\bibitem{cit:5} Coskun Hamzacebi, Fevzi Kutay, ``A heuristic approach for finding the global minimum: Adaptive random search technique," {\it Applied Mathematics and Computation}, vol. 173, pp. 1323-1333, 2006.
\bibitem{cit:6} Coskun Hamzacebi, Fevzi Kutay, ``Continous functions minimization by dynamic random search technique," {\it Applied Mathematical Modeling}, vol. 31, pp. 2189-2198, 2007.
\bibitem{cit:7} Tunchan Cura, ``A random search approach to finding the global minimum," {\it Int.J.Contemp.Math.Science}, vol. 5, no. 4, pp. 179-190, 2010.
\bibitem{cit:8} Xiaojun Zhou, ``Chunhua Yang and Weihua Gui, Initial version of State Transition Algorithm," {\it The 2nd International Conference on Digital Manufacturing \& Automation}, 2011(to be published).



\end{thebibliography}

\def\V{\rm vol.~}
\def\N{no.~}
\def\pp{pp.~}
\def\Pot{\it Proc. }
\def\IJCNN{\it International Joint Conference on Neural Networks\rm }
\def\ACC{\it American Control Conference\rm }
\def\SMC{\it IEEE Trans. Systems\rm , \it Man\rm , and \it Cybernetics\rm }

\def\handb{ \it Handbook of Intelligent Control: Neural\rm , \it
    Fuzzy\rm , \it and Adaptive Approaches \rm }

\end{document}